\documentclass[12pt]{article}
\usepackage{amsmath}
\usepackage{amscd}
\usepackage{graphicx}
\usepackage{amssymb,amsfonts,amsthm,amscd,mathrsfs}
\usepackage{verbatim}
\usepackage[T2A]{fontenc}
\usepackage[utf8]{inputenc}
\usepackage[russian]{babel}
\usepackage{bm}
\usepackage{authblk}

\renewcommand{\b}{\beta}

\newcommand{\e}{\varepsilon}

\newcommand{\E}{\mathsf{E}}

\begin{document}

\title{On Symmetrized Pearson's Type Test in Autoregression with Outliers: Robust Testing of Normality
}

\author{M.~V.~Boldin    
\footnote{Moscow State Lomonosov Univ., Dept. of Mech. and Math., Moscow, Russia\\
e-mail: boldin$_{-}$m@hotmail.com}}

\date{ }
\maketitle

\textbf{Abstract}

We consider a stationary linear AR($p$) model with observations subject to
gross errors  (outliers). The autoregression parameters are unknown as well as the distribution and moments of innoovations. 
The distribution of outliers $\Pi$ is unknown and arbitrary, their intensity is
 $\gamma n^{-1/2}$ with an unknown $\gamma$, $n$ is the sample size.
 The autoregression parameters are estimated by any estimator which is $n^{1/2}$-consistent
uniformly in $\gamma\leq \Gamma<\infty$. Using the residuals from
the estimated autoregression, we construct a kind of empirical distribution function
(e.d.f.), which is a counterpart of the (inaccessible)
e.d.f. of the autoregression innovations. We obtain a stochastic expansion
of this e.d.f., which enables us to construct the symmetrized test of Pearson's chi-square
type for the normality of distribution of innovations.
We establish qualitative robustness of these tests in terms of uniform equicontinuity
 of the limiting levels (as functions of $\gamma$ and $\Pi$) with respect to $\gamma$ in a neighborhood of $\gamma=0$.



\section{
Введение и постановка задачи}
В линейных и нелинейных регрессионных и авторегрессионных моделях при гауссовских ошибках (инновациях) обычно можно построить оптимальные или асимптотически оптимальные процедуры оценивания параметров и проверки гипотез относительно этих парамтров. Удается так же оптимально решать некоторые задачи типа change-point и разладки. Вот некоторые примеры. 

Классические результаты об оценивании и проверке гипотез в линейных AR, MA, ARMA моделях можно найти в \cite{And.}, \cite{Broc.Dav.}; результаты о GM-процедурах и процедурах минимального расстояния -- в \cite{Koul87} ; об оценивании в AR(1) с коэффициентом из $R^1$ и последовательных оценках максимального правдоподобия -- в \cite{Shir.Spok.}; результыты для задачи типа  change-point о дрейфе параметров в AR и ARCH моделях -- в \cite{Bold.Erl.}.

В связи со сказанным содержательной является задача проверки гипотезы о гауссовости инноваций в авторегрессионных моделях. Подобная и родственые задачи проверки гипотез о распределении инноваций критериями типа Колмогорова-Смирнова, Крамера-Мизеса-Смирнова омега-квадрат, хи-квадрат Пирсона  рассматривались для разных моделей, ссылки можно найти, нпример, в \cite{Bold.Petr.}, \cite{BST97}.

В этой работе мы рассматриваем стационарную AR($p$) модель с ненулевым средним
\begin{equation}
v_t = \b_1 v_{t-1} + \dots + \b_p v_{t-p} +\nu + {\e}_t, \quad  t \in \mathbb{Z,} 
\end{equation}

В (1.1) $\{{\e}_t\}$ -- независимые одинаково распределенные случайные величины (н.о.р.сл.в.) с неизвестной функцией распределения (ф.р.) $G(x)$;
 $\E \e_1 = 0$,  $0<\E {\e}_1^2 < \infty$; $\bm{\b} = (\b_1, \dots, \b_p)^T$-- вектор неизвестных параметров, таких что корни соответствующего (1.1) характеристического уравнения по модулю меньше единицы;  $\nu$ -- неизвестное среднее (константа из $R^1$).\\
 Эти условия дальше всегда предполагаютя выполненными и особо не оговариваются.\\
 В этой работе мы рассматриваем модель (1.1) с выбросами в наблюдениях. А именно, предполагается, что наблюдаются величины
$$
y_t=v_t+z^{\gamma_n}_t {\xi}_t,\quad t  = 1-p, \dots, n, \eqno(1.2)
$$
где $v_{1-p}, \dots, v_n$ -- выборка из стационарного решения $\{v_t\}$ уравнения  (1.1); $\{z^{\gamma_n}_t\}$ н.о.р.сл.в., принимающие значения 1 и 0, причем вероятность единицы
 $\gamma_n$,
 $$
 \gamma_n = \min(1, \frac{\gamma}{\sqrt{n}}) ,\quad  \gamma \ge 0 \text{\; неизвестно.}
 $$
Кроме того, $\{ \xi_t \}$ -- н.о.р.сл.в. с произвольным и неизвестным распределением $\Pi$.
Переменные $\{ \xi_t \}$ интерпретируются как выбросы (засорения), $\gamma_n$ уровень засорения. Для $\gamma = 0$ мы получаем модель (1.1) без засорений.

Модель (1.2) -- локальный вариант хорошо известной модели засорения данных во временных рядах, см. \cite{MartYoh86}.

Сформулируем цель настоящей работы.  Мы хотим проверить гипотезу о ноормальности инноваций

$$
H_{\Phi}\colon G(x) \in \{\Phi(x/\theta),\; \theta > 0\},
$$
где $\Phi(x)$ -- стандартная нормальная ф.р.

Напомним, кстати, что гипотеза $H_{\Phi}$ эквивалентна нормальности самой стационарной последовательности 
 $\{v_t\}$ itself.\\
 Подобная задача для авторегрессии с нулевым средним (т.е. схемы (1.1) с $\nu=0$) и засорениями в наблюдениях вида (1.2) рассматривалась в \cite{Bold.Petr.}. Для такой сиуации в  \cite{Bold.Petr.} по $\{y_t\}$ сначала строилась подходящая равномерно по $\gamma\leq \Gamma <\infty$ $n^{1/2}$-состоятельная оценка вектора $\bm{\b}$.
 Затем с ее помощью находились оценки ненаблюдаемых $\e_1,\ldots,\e_n$, они называются остатками. По ним строилась эмпирическая функция распределения остатков (о.э.ф.р.). Наконец, спомощью о.э.ф.р.  строилась подходящая статистика типа хи-квадрат Пирсона для $H_{\Phi}$. Эта статистика -- функционал о.э.ф.р. , это аналог обычной статистики Пирсона, построенной (гипотетически) по $\e_1,\ldots,\e_n$ и предназначенной для проверки нормальности при неизвестной дисперсии. Было найдено предельное распределение этой статистики при гипотезе. В схеме без засорений при $\gamma=0$ оно такое же, как у классической статистики Пирсона в теореме Фишера, построенной  по $\e_1,\ldots,\e_n$. При $\gamma\geq 0$ предельный уровень значимости соответствующего теста есть $\alpha(\gamma,\Pi),\,\alpha(0,\Pi)=\alpha$. Показано, что семейство $\{\alpha(\gamma,\Pi)\}$ равностепенно непрерывно по $\gamma$ в точке $\gamma=0$, т.е.
 $$
 \sup_{\Pi}|\alpha(\gamma,\Pi)-\alpha|\to 0,\quad \gamma\to 0.
 $$
 Последнее соотношение означает асимптотическую качествнную робастность теста относительно засорений. \\
 Недавно в  \cite{Bold.2020} была найдена асимптотическая локальная мощность теста из \cite{Bold.Petr.} и установлена ее качественная робастность.\\
 
 Упомянутые результаты из \cite{Bold.Petr.}-\cite{Bold.2020} получены при $\nu=0$ и это существенное предположение.
 В нашей ситуации  (1.1) при неизвестном и, быть может, ненулевом $\nu$, использовать статистику хи квадрат из \cite{Bold.Petr.}-\cite{Bold.2020} не удается, т.к. ее предельное распределение зависит от оценок $\bm{\b}$ и $\nu$ даже при $\gamma=0$ .\\
 Мы преодолеваем эту трудность, строя специальную симметризованную статистику Пирсона, она является функционалом от симметризованной эмпирической функции распределения остатков. Асимптотические распределения симметризованных статистик при $H_{\Phi}$ и  $\gamma=0$ свободны. Показано, что асимптотический уровень значимости нашего теста качественно робастен.
 
 Все определения и результаты (основные, это Теоремы 2.1--2.2) представлены в Разделе 2.
 
 \section{Основные результаты}

\subsection{Стохастическое разложение о.э. ф. р.}

Перепишем уравнение (1.1) в удобном для дальнейшего рассмотрения виде. Для этого определим константу $\mu$ соотношеним
$$
\nu=(1-\b_1-\ldots-\b_p)\mu,
$$
тогда
$$
v_t -\mu= \b_1 (v_{t-1}-\mu) + \dots + \b_p( v_{t-p} -\mu) {\e}_t, \quad  t \in \mathbb{Z.}
$$

Если положить $u_t:=v_t-\mu,$ то
$$
v_t=\mu+u_t,\quad u_t=\b_1 u_{t-1} + \dots + \b_pu_{t-p}+ {\e}_t, \quad  t \in \mathbb{Z.}\eqno(2.1	)
$$
Последовательность $\{u_t\}$ в (2.1) -- авторегрессионная последовательность с нулевым средним.\\
Построим по наблюдениям $\{y_t\}$ из (1.2) оценки ненаблюдаемых $\{\e_t\}$ в (2.1).\\ Пусть $\hat\mu_n$ будет любая последовательность, для которой последовательнсть
$$
n^{1/2}(\hat\mu_n-\mu)=O_P(1),n\to\infty, \text{ равномерно по}\,\,\gamma\leq \Gamma<\infty.
$$
 Положим
$$
\hat u_t=y_t-\hat\mu_n,\quad t=1,\ldots,n.
$$
Пусть $\hat{\bm\b}_n= (\hat \b_{1n}, \dots, \hat \b_{pn})^T$ будет любая последовательность, для которой последовательность\\
 $$
 n^{1/2}(\hat{\bm\b}_n-\bm\b)=O_P(1),n\to \infty, \text{ равномерно по}\,\,\gamma\leq \Gamma.
 $$
Примеров подходящих оценок $\hat\mu_n,\,\hat{\bm\b}_n$ много, широкий класс составляют, например, М-оценки. Результаты для них аналогичны изложенным в Разделе 2.4 \cite{Bold.Petr.}.\\
Положим
$$
\hat \e_t =  \hat {u}_t -  \hat \b_{1n}\hat { u}_{t-1} - \dots - \hat \b_{pn}\hat{ u}_{t-p},\quad t = 1, \dots, n,
$$
Наша ближайшая цель -- исследовать асимптотические свойства остаточной эмпирической ф.р.
$$
\hat G_n(x) = n^{-1} \sum_{t=1}^n I(\hat \e_t \le x),\quad  x\in \mathbb{R}^1,
$$
 Здесь и в дальнейшем $I(\cdot)$ обозначает индикатор события.\\
Функция $\hat G_n(x)$ -- аналог гипотетической э.ф.р.
$$
G_n(x) = n^{-1} \sum_{t=1}^n I(\e_t \le x)
$$
ненаблюдаемых величин $\e_{1}, \dots, \e_n$.

\newtheorem{Th}{Теорема}[section]
\begin{Th}
Предположим, что ф.р. $G(x)$ дифференцируема с производной (плотностью) $g(x)$. Предположим, что $\sup_x|g'(x)| < \infty$. Пусть $\Gamma \ge 0$ and $\Theta \ge 0$
будут любое действительное число. Тогда для любого $\delta > 0$
$$
\sup_{\gamma \le \Gamma} \Prob(|n^{1/2} [\hat G_n(x) - G_n(x)] -g(x)(1-\b_1-\ldots-\b_p) n^{1/2}(\hat\mu_n-\mu)- \gamma\Delta(x, \Pi) | > \delta) \to 0,\quad n \to \infty.
$$
Здесь сдвиг
$$
\Delta(x, \Pi) = \sum_{j=0}^p [\E G(x + \b_j \xi_j) - G(x)],\, \b_0 =-1.
$$
\end{Th}
Доказательство Теоремы 2.1 проводится по той же схеме, что доказательство Теореы 2.1 в \cite{Bold.Petr.}.

В связи с проверкой гипотезы $H_{\Phi}$ нас особо интересует ситуация, когда ф.р. $G(x)$ симметрична относительно нуля. В этом случае будем брать оценкой $G(x)$ симметризованную оценку
$$
\hat S_n(x):=\frac{\hat G_n(x)+1-\hat G_n(-x)}{2}.
$$
Положим
$$
\Delta_S(x, \Pi):=\frac{\Delta(x, \Pi)-\Delta(-x, \Pi)}{2}.\eqno(2.2)
$$
Пусть
$$
 S_n(x):=\frac{ G_n(x)+1- G_n(-x)}{2}.
$$

Теорема 2.1 влечет
\newtheorem{Corollary}{Следствие
}[section]
\begin{Corollary}
При условиях Теоремы 2.1 в случае симметричной относительно нуля $G(x)$
$$
\sup_{\gamma \le \Gamma} \Prob(|n^{1/2}[\hat S_n(x) -S_n(x)]-  \gamma\Delta_S(x, \Pi) | > \delta) \to 0,
\quad n \to \infty.
$$
\end{Corollary}

\subsection{Тест типа хи-квадрат Пирсона для $H_{\Phi}$}

В этом разделе мы построим тест типа хи-квадрат Пирсона для проверки 
$$
H_{\Phi}\colon  G(x) \in \{\Phi(x/\theta),\;\theta>0\},
$$
где $\Phi(x)$ стандартная нормальная ф.р. Как отмечалось во Введении, такая гипотеза эквивалентна гауссовости стационарного решения (1.1).

Обозначим через $\theta_0$ истинное значение $\theta$ при $H_{\Phi}$, тогда при $H_{\Phi}$ $G(x)=\Phi(x/\theta_0)$.
Для полуинтервалов 
$$
B^+_j=(x_{j-1},x_j],\quad j=1,\ldots,m,\; m>2,\; 0=x_0<x_1<\ldots<x_m=\infty,
$$
пусть 
$$
p^+_j(\theta)=\Phi(x_{j}/\theta)-\Phi(x_{j-1}/\theta).
$$
При $H_{\Phi}$ $\Prob(\e_1 \in B^+_j)=p^+_j(\theta_0)$. Если ввести еще симметричные полуинтервалы

 $ B^-_j=(-x_{j},-x_{j-1}]$, то при $H_{\Phi}$
$$
\Prob(\e_1 \in B^+_j \cup B^-_j)=2p^+_j(\theta_0):=p_j(\theta_0).
$$
Пусть $\hat \nu_j^+$ обозначает число остатков среди $\{\hat \e_t,\,t=1,\ldots,n\}$, попавших в $B^+_j$, а 
$\hat \nu_j^-$ обозначает число остатков, попавших в $B^-_j$. Пусть
$$
\hat\nu_j=\hat \nu_j^+ +\hat \nu_j^-.
$$
Оценкой $\theta_0$ мы возьмем $n^{1/2}$-состоятельное решение уравнения
$$
\sum_{j=1}^m \frac{\hat{\nu}_j}{p_j(\theta)}p_j^{'}(\theta)=0, \eqno(2.3)
$$
в котором штрих означает производную относительно  $\theta$. Такое уравнение -- аналог обычного уравнения модифицированного метода хи-квадрат, в котором ненаблюдаемые частоты $\nu_1,\ldots,\nu_m$ инноваций, попавших в $B_1,\ldots,B_m$, заменены их оценками $\hat{\nu}_1,\ldots,\hat{\nu}_m$.

Интересующая нас тестовая статистика типа хи-квадрат для $H_{\Phi}$ имеет вид
$$
\hat{\chi}^2_{n} = \sum_{j=1}^m \frac{(\hat{\nu}_j - n p_j(\hat{\theta}_n))^2}{n p_j(\hat{\theta}_n)}.\eqno(2.5)
$$
Статистика $\hat{\chi}^2_{n} $ является простым функционалом от $\hat S_n(x)$. Действительно,

$$
n[\hat S_n(x_j)-\hat S_n(x_{j-1})]=\frac{1}{2}\{n[\hat G_n(x_j)-\hat G_n(x_{j-1})]+n[\hat G_n(-x_{j-1})-\hat G_n(-x_{j})]\}=
$$
$$
(\hat\nu^+_j +\hat\nu^-_j)/2)=\hat\nu_j/2.\eqno(2.4)
$$
Будем называть $\hat{\chi}^2_{n}$ симметризованной статистикой Пирсона.\\
Чтобы описать асимптотические свойства $\hat{\chi}^2_{n}$ нам понадобятся некоторые обозначения.
Введем диагональную матрицу
$$
\bm P(\theta)=\mbox{diag}\{p_1(\theta),\ldots,p_m(\theta)\}
$$
и вектора
$$
\bm p(\theta)=(p_1(\theta),\ldots,p_m(\theta))^T,\quad
\bm p'(\theta)=(p_1'(\theta),\ldots,p_m'(\theta))^T,
$$
$$
\bm b(\theta)=\bm P^{-1/2}(\theta)\bm p'(\theta), \quad \bm \alpha(\theta)=\bm b(\theta)/|\bm b(\theta)|.
$$
Здесь и далее $|.|$ означает евклидову норму вектора или матрицы.\\
Нам понадобится еще вектор
$$
\bm \delta(\Pi)=(\delta_1(\Pi),\ldots,\delta_m(\Pi))^T\quad \text{с компонентами}\,\, \delta_j(\Pi)=2[\Delta_{S}(x_{j},\Pi)-\Delta_{S}(x_{j-1},\Pi)],
$$
$\Delta_{S}(x,\Pi)$ определена в (2.2), где следует положить  $G(x)=\Phi(x/\theta_0)$.

Из определения $S_n(x)$ и $\nu_j$ следует:
$$
n[ S_n(x_j)- S_n(x_{j-1})]=\nu_j/2.\eqno(2.5)
$$
В силу соотношений (2.4) --(2.5), Следствия 2.1 и определения вектора $\delta(\Pi)$ имеем при $H_{\Phi}$:
$$
n^{1/2}(\frac{\hat\nu_j}{n}-
p_j(\theta_0))-n^{1/2}(\frac{\nu_j}{n}-p_j(\theta_0))=
n^{1/2}(\frac{\hat\nu_j}{n}-\frac{\nu_j}{n})=
$$
$$
2n^{1/2}\{[\hat S_n(x_j)-\hat S_n(x_{j-1}]-[S_n(x_j)-S_n(x_{j-1}]\}=\gamma \delta_j(\Pi)+o_P(1),\quad n\to \infty.\eqno(2.6)
$$
Здесь $o_P(1)$ обозначает величину, сходящуюся к нулю по вероятности равномерно по $\gamma \leq \Gamma$.\\
Введем вектора 
$$
\bm\hat \nu=(\hat\nu_1,\ldots,\hat\nu_n)^T,\quad \bm \nu=(\nu_1,\ldots,\nu_n)^T.
$$
В силу (2.6) при $H_{\Phi}$
$$
n^{1/2}(\frac{\bm \hat \nu}{n}-\bm p(\theta_0))=
n^{1/2}(\frac{\bm \nu}{n}-\bm p(\theta_0))+\gamma \bm\delta(\Pi)+o_P(1),\quad n\to \infty.\eqno(2.7)
$$
Напомним еще общеизвестный факт относительно слабой сходимости вектора $\bm \nu$:
$$
n^{1/2}(\frac{\bm \nu}{n}-\bm p(\theta_0)) \to \bm N(\bm 0,\,\, \bm P(\theta_0)-\bm p(\theta_0)\bm p(\theta_0)^T),\quad n\to \infty.\eqno(2.8)
$$
Соотношений (2.7) --(2.8) достаточно, чтобы следуя схеме доказательства Теоремы 2.3 в \cite{Bold.Petr.} доказать наше  основное утверждение -- Теорему 2.2. В ней $\bm E_m$ означает единичную матрицу порядка  $m$.
\begin{Th}
Пусть $H_{\Phi}$ верна. Тогда для любого конечного  $0 \le \Gamma < \infty$
$$
\sup_{x\in\mathbb{ R}^1, \gamma \le \Gamma} |\Prob(\hat{\chi}^2_{n} \le x) - F_{m-2}(x, \hat\lambda^2(\gamma, \Pi))| \to 0,\quad n \to \infty.
$$
Параметр нецентральности равен
$$
\hat\lambda^2(\gamma, \Pi) = \gamma^2
|(\bm{ E}_m- \bm \alpha(\theta_0) \bm\alpha(\theta_0)^T)
 \bm{ P}^{-1/2}(\theta_0) \bm{\delta}(\Pi)|^2.
$$
\end{Th}
Мы будем отвергать $H_{\Phi}$, когда
$$
\hat{\chi}^2_{n} > \chi_{m-2}(1-\alpha), \eqno(2.9)
$$
$\chi_{m-2}(1-\alpha)$ -- квантиль уровня $1-\alpha$ распределения хи-квадрат с $m-2$ степенями свободы.\\
В силу Теоремы 2.2 уровень такого теста сходится при $n\to\infty$ равномерно по $\gamma\leq \Gamma$
к асимптотическому уровню
$$
\hat{\alpha}( \gamma,\Pi) = 1 - F_{m-2}(\chi_{m-2}(1-\alpha), \hat\lambda^2(\gamma, \Pi)),\quad
\hat \alpha(0,\Pi) = \alpha.
 $$
 Легко проверить, что
 $$
\sup_{\Pi}|\hat { \alpha}( \gamma,\Pi) -\alpha| \to 0,\quad \gamma \to 0.\eqno(2.10)
$$
Соотношение (2.10) означает асимптотическую качественную робастность теста (2.9). Такая робастность означает, что  
 $H_{\Phi}$ для малых  $\gamma$ можно проверять примерно с асимптотическим уровнем $\alpha$ независимо от распределения засорений $\Pi$.\\

\end{document}